\documentclass{article}

\usepackage[margin=1.5in]{geometry}
\usepackage[utf8]{inputenc} 
\usepackage[T1]{fontenc}    
\usepackage{hyperref}       
\usepackage{url}            
\usepackage{booktabs}       
\usepackage{amsfonts}       
\usepackage{nicefrac}       
\usepackage{microtype}      
\usepackage{graphics} 
\usepackage{amsthm}
\usepackage{epsfig} 
\usepackage{tikz,cite}
\usepackage{amsmath} 
\usepackage{amssymb}  
\usepackage{setspace}
\usepackage[version=4]{mhchem}
\graphicspath{{./Figures/}}

\newtheorem{theorem}{Theorem}

\usepackage{url}

\title{A Data-Driven Sparse-Learning Approach to Model Reduction in Chemical Reaction Networks}

\author{ 
  Farshad Harirchi$^a$, Omar A. Khalil$^a$, Sijia Liu$^a$, Paolo Elvati$^c$,\\ Angela Violi$^{c,d}$,Alfred O. Hero$^{a,b}$\\ 
  $^a$Department of Electrical Engineering and Computer Science,\\ $^b$Department of Statistics,\\
  $^c$Department of Mechanical Engineering,\\ 
  $^d$Departments of Chemical and Biomedical Engineering,\\
   University of Michigan, Ann Arbor, MI 48109 \thanks{This research has been funded by the US Army Research Office grants $\#$W911NF-15-1-0241 and $\#$W911NF-14-1-0359.}\\
  \texttt{$\{$harirchi,oakhalil,lsjxjtu,elvati,avioli,hero$\}$@umich.edu} \\
}

\begin{document}

\maketitle

\begin{abstract}
In this paper, we propose an optimization-based {\em sparse learning} approach to identify the set of most influential reactions in a chemical reaction network. This reduced set of reactions is then employed to construct a reduced chemical reaction mechanism, which is relevant to chemical interaction network  modeling. The problem of identifying influential reactions is first formulated as a mixed-integer quadratic program, and then a relaxation method is leveraged to reduce the computational complexity of our approach. Qualitative and quantitative validation of the sparse encoding approach demonstrates that the model captures important network structural properties with moderate computational load.
\end{abstract}

\section{Introduction} 
Many industrial and natural processes are the result of thousands of
interdependent chemical reactions that happen on very short time
(femto- to millisecond) scales.
Despite the increase in computational power of the last few decades,
the modeling of these Chemical Reaction Networks (CRN) for realistic
applications remains a daunting task, as the number of reactions in a mechanism can exceed several tens of thousands.
The simultaneous resolution of such number of differential equation is only possible under very simple conditions, like zero-dimensional systems, but not in general, specifically for applications such as fluidodynamics.
To overcome this problem, various reduced versions of the chemical
kinetic mechanisms have been developed [3-7].

The goal of these reduction methods is to identify a small subset of species
and reactions that are able to approximate complete CRN with a
drastic reduction in the computational complexity.
As the evolution of chemical reactions depends on a set of boundary
conditions (such as temperature, pressure, and equivalence ratio) [2], the topology of CRN changes in time. This paper adopts a data-driven sparse learning approach to model reduction of CRN.

\textbf{Related Work:} The problem of mechanism reduction in chemical reaction networks has been considered in the combustion community. In [3] a method based on Jacobian analysis is proposed, which identifies species that are coupled with the important species. This method requires a good knowledge of the reaction network, choice of important species, and arbitrary choice of threshold values. Furthermore, it is an iterative approach that is computationally demanding. A commonly used approach for mechanism reduction is developed by utilizing the Directed Relational Graphs (DRG) to first identify the skeletal mechanism [4,5] based on a normalized measure of direct influence of one species on the production rate of another. It utilizes a threshold-based technique to add a directed edge from one species to another. In the second step, the mechanism is further reduced by using a quasi steady state assumption (QSS) [6]. This method is strongly mechanism dependent, and requires a very good understanding of the process. In addition, it can be inefficient to remove modes at short time scales. On the other hand, the choice of threshold on each step is arbitrary. Finally, in [7] a revised version of DRG approach is proposed that takes into account the error propagation effect by revising the normalized measure of direct influence.

\textbf{Contributions:} In this paper, we propose an optimization-based sparse learning approach to the problem of selecting the most influential reactions. First, a data-driven distance measure is defined between the detailed mechanism and the reduced mechanism. The learner finds the minimal set of reactions that keeps the fitting error, i.e., the distance between the two mechanisms, within some range for all time instances for all the data points. Then, we calculate the set union of all the selected reactions at all time instances, and for a sparse set of points in the parameter space. Unlike other methods in the literature, the proposed method does not require any understanding of the underlying chemical process, instead it learns the model from data generated from the chemical process. Additionally, our method only needs tuning of a single parameter, which is the tolerance on the distance. We verify the efficacy of our approach on the $\ce{H2}/\ce{O2}$ reaction network by verifying the reduced mechanism on the data obtained from the remaining points in the parameter space. We propose an $\ell_1$-sparse learning approach to relax the mixed integer quadratic program to a quadratic program that can be efficiently solved on a much larger scale.

\section{Modeling Framework} 
The dynamics of chemical reaction networks, under some conditions, can be reduced to a set of ordinary differential equations that are known as the mass-action kinetic equations [1]. For a mechanism with $N_s$ species and $N_r$ reactions, these equations are nonlinear ODE's of the following form:  \vspace*{-0.2cm}
\begin{equation}\vspace*{-0.2cm}
\frac{d}{d t} \mathbf{X}_t = \mathbf{M} \mathbf{r}_t, \; t \in \mathcal{\tau}, \mathbf r_t = f(\mathbf{X}_t)
\end{equation}  
where $\mathbf{X}_t\in \mathbb{R}_{\geq0}^{N_s}$ denotes the vector of molar concentrations of all species at time $t$, $\mathbf{r}_t \in \mathbb{R}_{\geq0}^{N_r}$ represents the reaction rates for all reactions at time $t$, $\mathbf{M} \in \mathbb{Z}^{N_s \times N_r}$ contains the stoichiometric coefficients corresponding to all reactions, $\mathcal{\tau}$ is the time interval of interest and $f(.)$ is a polynomial function operating on individual species concentrations as described in [1]. Our approach is agnostic to the specific form of the nonlinearity introduced by $f(.)$ since it works directed on data generated by a {\em generative } model. By employing the Euler discretization method, we can discretize the continuous dynamics with the following model:  \vspace*{-0.2cm}
\begin{equation} \label{eqn:discrete_model} \vspace*{-0.2cm}
\mathbf{X}_{t+1} = \mathbf{X}_t + \sum_{i=1}^{N_r} \mathbf{M}_i \mathbf{r}_{i,t} \Delta t + \pmb{w}_t, \; t = 1,2, \hdots , T,
\end{equation}  
where $\Delta t$ is the sampling time, $N_r$ is the number of reactions, and $\pmb{w}_t \in \mathbb{R}^{N_s}$ denotes the discretization error and process noise. 

\section{Sparse Learning Approach for Mechanism Reduction} 
In this section, we propose a data-driven sparse learning approach to identify a reduced set of reactions, which approximately replicates the behavior of the detailed mechanism. Our learner finds the smallest set of reactions that keeps a measure of distance (the fitting error) between the reduced and detailed mechanisms within some error bound $\epsilon$. The reduced mechanism, then only depends on the choice of $\epsilon$. A larger value for $\epsilon$ results in a more reduced mechanism at the cost of larger distance. 
The fitting error is defined as the Euclidean
distance for the reduced mechanism $\mathcal{G}$ at time $t$:\vspace*{-0.1cm}
\begin{equation} \label{eqn:distance}\vspace*{-0.1cm}
\mathcal{D}_t(\mathcal{G}) = \|\mathbf X_{t+1}-\mathbf X_t - \mathbf M \mathbf r^\mathcal{G}_t \Delta_t\|_2, \quad \mathbf r^{\mathcal{G}}_t = \mathbf w_t \odot \mathbf r_t,
\end{equation} 
where $\mathcal G$ represents a reduced mechanism parameterized by a binary selection vector $\mathbf w_t$, and $\odot$ denotes the element-wise product. Note that if there is no model mismatch or process noise, then $\mathcal{D}_t(\mathcal{G})$ denotes the Euclidean distance between the original and reduced mechanisms $\mathcal{G}$ at time $t$. The fitting error measure $\mathcal{D}_t(\mathcal{G})$ in \eqref{eqn:distance} can also be interpreted as a denoising function under the linear model \eqref{eqn:discrete_model}. Here the variable $\mathbf w_t$ is introduced to encode whether or not a reaction is selected, e.g., $\mathbf w^i_t=1$  if reaction $i$ is selected and $0$ otherwise, where $\mathbf w^i_t$ denotes the $i^{th}$ entry of $\mathbf w_t$.
We can control fitting errors by imposing the following constraint that determines the goodness of fit of the reduced set of reactions.  \vspace*{-0.2cm}
\begin{equation} \label{eqn:distancecriteria}\vspace*{-0.1cm}
\mathcal{D}_{t}(\mathcal{G}(\mathbf w_t))\leq \frac{\epsilon}{\mathcal{N}_t}  , \forall t\in \{1,2, \hdots , T \},\; \mathcal{N}_t=\|\mathbf 1 \mathbf r_t \Delta_t\|_2.
\end{equation} 
Equation \eqref{eqn:distancecriteria} ensures that the fitting error does not exceed $\epsilon$ normalized by the $\ell_2$-norm of reaction rates at time $t$, namely $\mathcal{N}_t$, at any given time. In other words, the error in the change of concentration introduced by reducing the mechanism remains in some percentage ($\epsilon$) of the norm of reaction rate at each time. 
In the presence of model mismatch or process noise, 
equation \eqref{eqn:distancecriteria} can be interpreted as a constraint that limits the robustness to modeling error and process noise similar to robust compressive sensing problems [8].

The objective of our approach is to identify the minimal set of reactions that satisfy the distance criterion at all time instances. To achieve this goal, we formulate the following optimization problem. 
\begin{align}\label{eq:prob_ori}
\begin{array}{ll}
\min_{\mathbf{w}_t} & \sum_{i=1}^{N_r} \mathbf{w}_t^i\\
\mathrm{subject} \; \mathrm{to} & \mathcal{D}_{t}(\mathcal{G}(\mathbf w_t))\leq \frac{\epsilon}{\mathcal{N}_t},\\
& \mathbf w_t^i \in \{0,1\}, ~i=1,2,...,N_r,
\end{array}
\end{align}
where $\mathbf w_t$ is a binary selection vector, as defined above. Even though, the worst-case complexity of the mixed-integer quadratic programming problem defined by \eqref{eq:prob_ori} is exponential in the integer variables, this problem can be solved efficiently using state of the art MIQP solvers such as Gurobi [9]. After finding the minimal set of influential reactions at each time instance and for a few (say $J$) points in the parameter space, we find the set union of all these reactions over the whole time horizon and all selected points in the parameter space as follows:
\begin{equation} \label{eqn:union}
\begin{aligned}
\mathbf w = \cup_{j=1}^{J}\cup_{t=1}^{T} \mathcal{R}(\mathbf w_{t,j}), \; \mathcal{R}(\mathbf w_{t,j}): & \text{ specified by selection vector $\mathbf w_{t,j}$}.
\end{aligned}
\end{equation}
This gives us the minimal set of reactions required to satisfy the distance criterion at all times and for all the points selected in the parameter space. 

\begin{theorem}
	The reduced mechanism with the set of reactions obtained from \eqref{eqn:union} satisfies the fitting error criterion.
\end{theorem} \vspace*{-0.5cm}	
\begin{proof}
	Note that from \eqref{eqn:union}, the following element-wise inequality holds:
	\begin{equation}
	\mathbf w^i \geq \mathbf w_{t,j}^i, \forall t,\forall i\in \{1, \hdots ,N_r\}, \forall j\in \{1, \hdots ,J\}
	\end{equation} 
	therefore, from \eqref{eqn:distance}, we always have $\mathcal{D}_t(\mathcal{G(\mathbf w)}) \leq \mathcal{D}_t(\mathcal{G}(\mathbf w_{t,j}))$ at all time instances and for all selected conditions on the parameter space. This concludes that
		
		$\qquad \qquad \qquad \qquad \quad \mathcal{D}_t(\mathcal{G(\mathbf w)}) \leq \frac{\epsilon}{\mathcal{N}_t}, \; \forall t,\forall i\in \{1, \hdots ,N_r\}, \forall j\in \{1, \hdots ,J\}.$\centering
		\end{proof}
\subsection{Relaxed Mechanism Reduction Approach} 
The proposed mechanism reduction approach involves solving $J\times T$ mixed integer quadratic programming problems, each of which has $N_r$ binary variables. Even though solving each of these problems can be performed in a tractable manner, if the size of mechanism grows, then this approach requires a significant amount of resources. As such, we consider a relaxed version of this problem that reduces the computational complexity and can be performed efficiently for large-scale problems. The relaxed problem has the following form:
\begin{align}\label{eq:prob_relaxed}
\begin{array}{ll}
\min_{\mathbf{w}_t} & \sum_{i=1}^{N_r} \mathbf{w}_t^i\\
\mathrm{subject} \; \mathrm{to} & \mathcal{D}_{t}(\mathcal{G}(\mathbf w_t))\leq \frac{\epsilon}{\mathcal{N}_t},\\
& 0 \leq \mathbf w_t^i \leq 1, ~i=1,2,...,N_r,
\end{array}
\end{align}
The real valued $\mathbf{w}_t$ obtained as the solution of problem \eqref{eq:prob_relaxed} is not necessarily feasible for the original problem \eqref{eq:prob_ori}, however, it can be interpreted as the normalized influence of each reaction at time $t$. By employing thresholding or randomization methods [10-11], we can readily project this solution to the feasible set of the original problem. 
To further reduce the mechanism, we can omit the reactions that only appear in a few time instances. Removing these reaction will result in the violation of distance criterion for those few time instances, but it may significantly reduce the size of mechanism.

\vspace*{-0.3cm}
\begin{figure}[!h]
	\centering
	\includegraphics[width=1.8in]{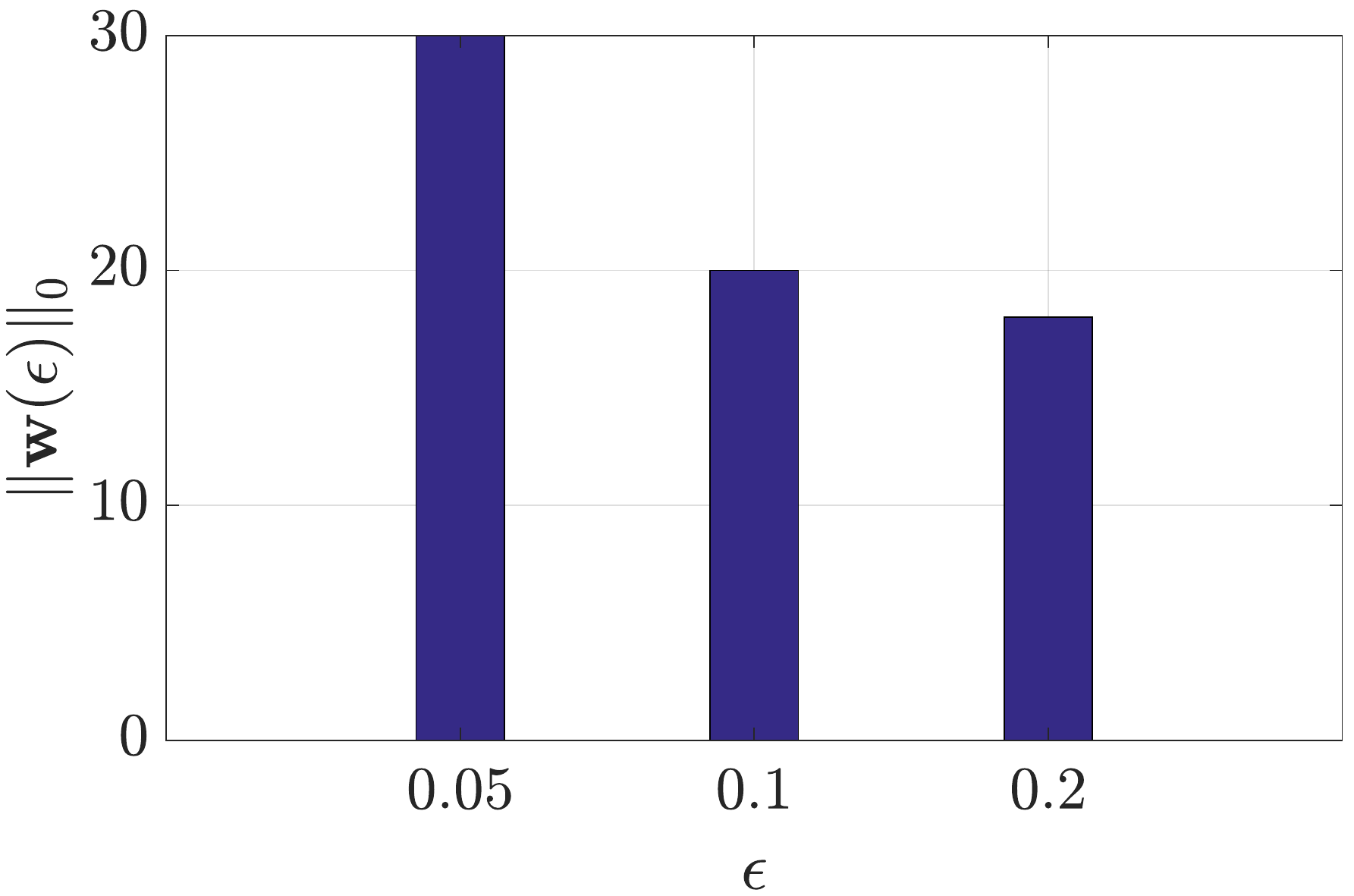}
	\includegraphics[width=1.8in]{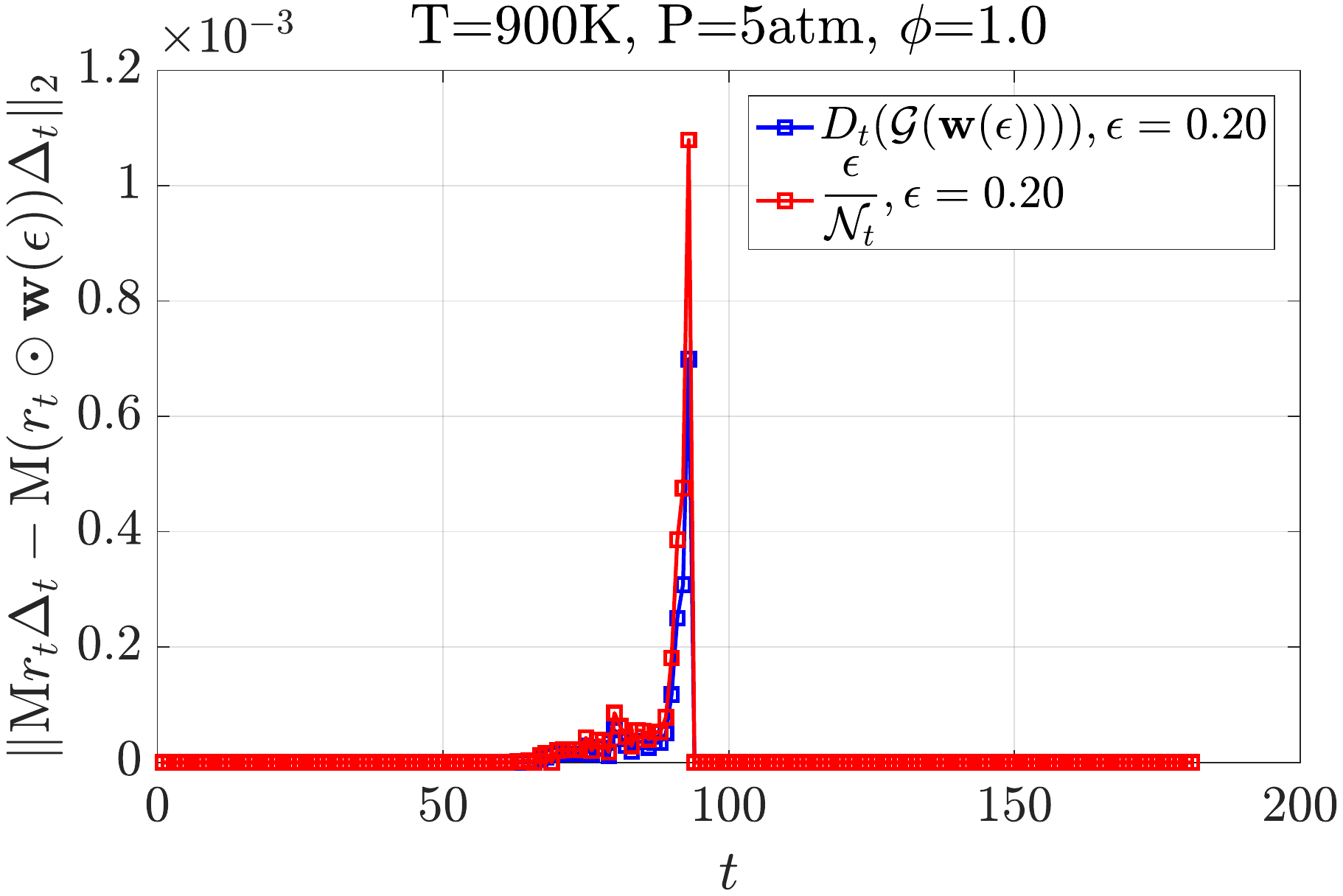}\\
	\includegraphics[width=1.8in]{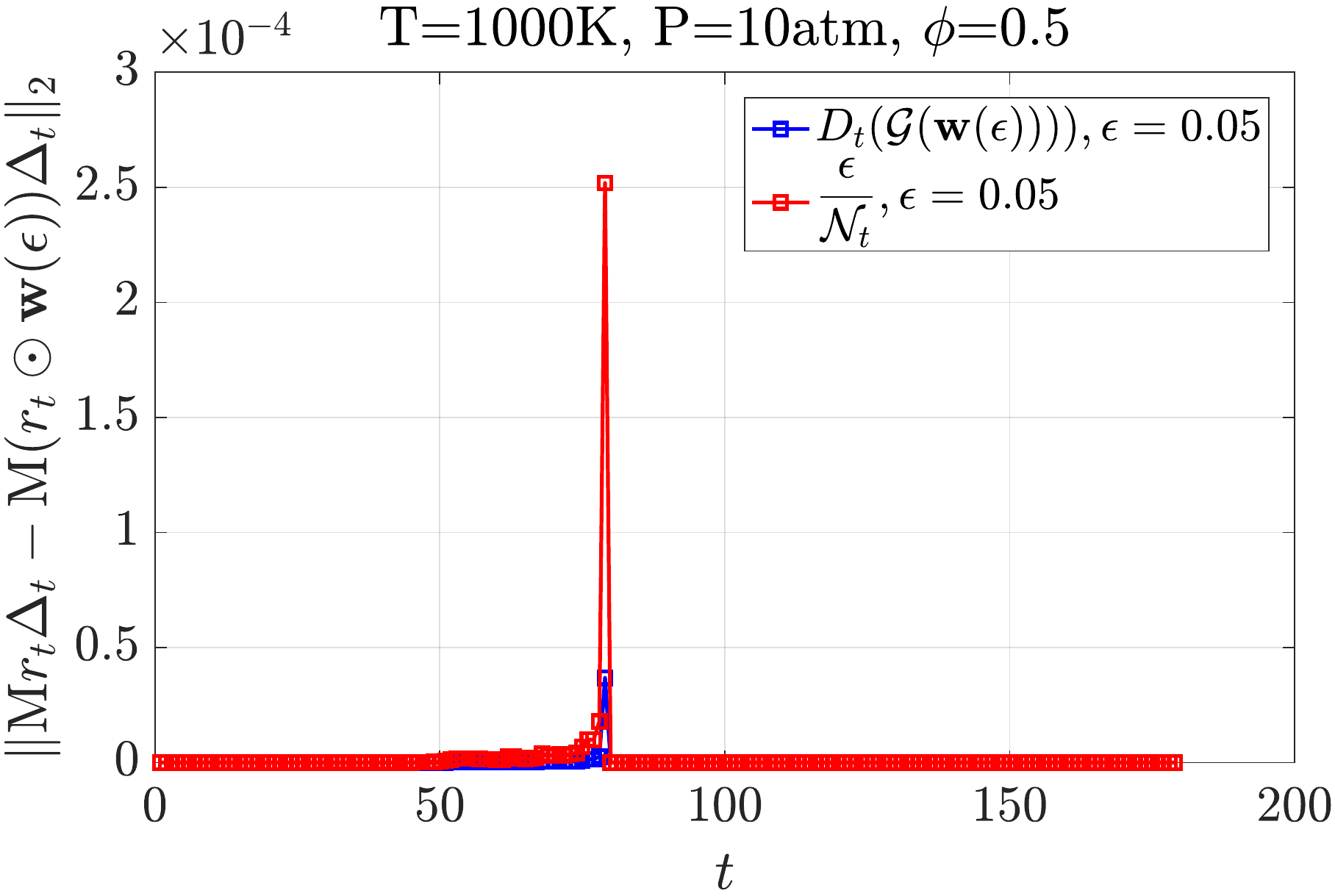} 
	\includegraphics[width=1.8in]{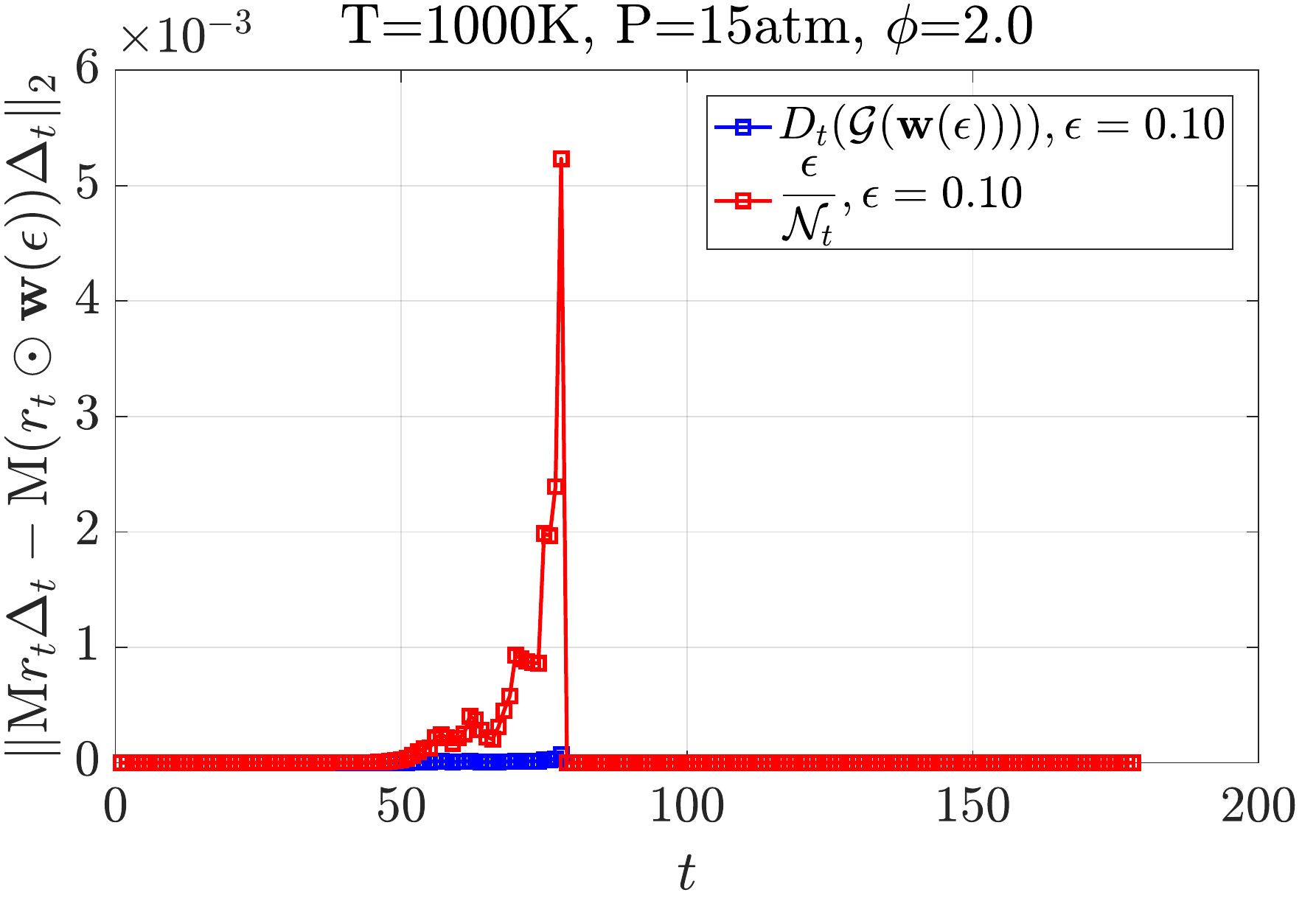} \vspace*{-0.3cm}
	\caption{Top Left: shows that with increase in error tolerance parameter $\epsilon$, the complexity of the mechanism is reduced. Top Right: illustrates the fitting error for a condition that is used to construct the reduced mechanism for $\epsilon = 0.2$ and its upper-bound, which is relatively tight. Bottom: the fitting error (blue curve) is compared to the upper-bound (red curve) for the two initial conditions that are not in the training set used for calculation of the reduced mechanism. As shown, the fitting error of the reduced mechanism resides within the bounds even for these initial conditions that are outside of the training set.} \label{fig:reduced}
\end{figure}

\section{Simulation Results} 
In this section, we illustrate the efficacy of the proposed mechanism reduction approach on a dataset generated by CHEMKIN\footnote{\url{http://www.reactiondesign.com/products/chemkin}} for 48 different initial conditions (i.e., temperature $(T)$, pressure $(P)$ and equivalence ratio $(\phi)$) on the \ce{H2}/\ce{O2} reaction mechanism that involves $8$ species and $58$ reactions. We implemented the sparse learning approach (i.e., \eqref{eq:prob_ori}) in Matlab using YALMIP [12] with Gurobi [9]. We selected three different values for a design parameter $\epsilon$, namely $0.05, 0.1$ and $0.2$. The number of selected reactions over $46$ initial conditions used for training and the fitting error for one initial condition are illustrated in Fig. \ref{fig:reduced}.
In order to validate the performance of our reduced mechanism, we further illustrate the fitting error and the corresponding upper bound for two different values of $\epsilon$ for the two initial conditions that are not contained in the training set.

\section{Conclusions} 
In this paper, we proposed a sparse learning optimization-based framework to the problem of mechanism reduction in chemical reaction networks. Our approach does not require full understanding of the mechanism, and is tuned with only one parameter that reflects the fitting error tolerance bound. Because of the data-driven nature of this approach, it can be implemented on any reaction mechanism. In order to reduce the computational complexity of our approach, we formulated a relaxed version of the optimization problem that is scalable for dealing with large CRN. To further improve the proposed sparse-learning approach, we can also constraint the propagation of the fitting error over some  time horizon. In addition, by employing  $\ell_{\infty}$ distance instead of $\ell_{2}$, we can replace the quadratic constraints with linear ones and therefore, reduce the complexity of the method. 

\section*{References}
\small
[1] Feinberg, M., (1987) Chemical reaction network structure and the stability of complex isothermal reactors–I: The deficiency
zero and deficiency one theorems. {\it Chemical Engineering Science 42}, pp. 2229--2268.

\vspace*{-0.11cm}
[2] Strehlow, R.A., (1984) Combustion fundamentals, {\it McGraw-Hill College}.

\vspace*{-0.11cm}
[3] Tomlin, A.S. \ \& Turanyi, T. \ \&  Pilling, M.J., (1997) Mathematical Tools for the Construction, Investigation and Reduction of Combustion Mechanisms. {\it Elsevier Comprehensive Chemical Kinetics}, pp. 293--437.

\vspace*{-0.11cm}
[4] Lu, T.\ \& Law, C.K., (2005) A directed relation graph method for mechanism reduction. {\it Elsevier Proceedings of the Combustion Institute 30}, pp. 1333--1341.

\vspace*{-0.11cm}
[5] Lu, T.\ \& Law, C.K., (2006) On the applicability of directed relation graphs to the reduction of reaction mechanisms. {\it Elsevier Combustion and Flame 146}, pp.\ 472--483.

\vspace*{-0.11cm}
[6] Tomlin, A.S.\ \& Pilling,  M.J. \ \& Turanyi, T.\ \& Merkin, J.H.\ \& Brindley, J., (1992). Mechanism Reduction for the Oscillatory Oxidation of Hydrogen: Sensitivity and Quasi-Steady-State Analyses. {\it Elsevier Combustion and Flame 91}, pp. 107--130.

\vspace*{-0.11cm}
[7] Pepiot-Desjardins, P.\ \& Pitsch, H., (2008) An efficient error-propagation-based reduction method for large chemical kinetic mechanisms. {\it Elsevier Combustion and Flame 154}, pp.\ 67--81.

\vspace*{-0.11cm}
[8] Candes, E.J.\ \& Wakin, M.B., (2010). An Introduction To Compressive Sampling. {\it IEEE Signal Processing Magazine}, pp. 21--30.	

\vspace*{-0.11cm}
[9] Inc. Gurobi Optimization, (2015) Gurobi optimizer reference manual.

\vspace*{-0.11cm}
[10] Joshi, S.\ \& Boyd, S., (2009) Sensor selection via convex optimization. {\it IEEE Transactions on Signal Processing 57}. pp. 451--462.

\vspace*{-0.11cm}
[11] Liu, S.\ \& Chepuri, S.P.\ \& Fardad M., (2016) Sensor selection for estimation with correlated measurement noise. {\it IEEE Transactions on Signal Processing 64}. pp. 3509--3522.
	
\vspace*{-0.11cm}
[12] Löfberg J., (2004). Yalmip : A toolbox for modeling and optimization in MATLAB. {\it In CACSD}.

%
%
%
%

\end{document}